\documentclass{amsart}
\usepackage{amsmath,amsfonts}

\usepackage{pstricks,pst-node,pst-plot}
\usepackage{epic}
\usepackage{autograph}

\newtheorem{theorem}{Theorem}
\newtheorem{proposition}[theorem]{Proposition}

\newtheorem{question}[theorem]{Question}
\newtheorem{lemma}[theorem]{Lemma}

\newtheorem{corollary}[theorem]{Corollary}
\newcommand{\ol}[1]{\overline{#1}}

\newcommand{\Z}{\mathbb{Z}}

\begin{document}

    \title[Graph groups and monoids]{On commuting elements and embeddings of\\ graph groups and monoids}

\maketitle

\begin{center}
 \maketitle

    Mark Kambites \\

    \medskip

    Fachbereich Mathematik / Informatik, \  Universit\"at Kassel \\
    34109 Kassel, \  Germany

    \medskip

    \texttt{kambites@theory.informatik.uni-kassel.de} \\

\end{center}

\begin{abstract}
We study commutation properties of subsets of right-angled Artin groups
and trace monoids. We show that if $\Gamma$ is any graph not containing
a four-cycle without chords, then the group $G(\Gamma)$ does not contain
four elements whose commutation graph is a four-cycle; a consequence is
that $G(\Gamma)$ does not have a subgroup isomorphic to a direct product
of non-abelian groups. We also obtain corresponding and more general
results in the monoid case.
\end{abstract}

\section{Introduction}\label{sec_intro}

Much research has centred upon finitely generated monoids
and groups defined by presentations in which the only relations are
commutators of certain of the generators.
Monoids of this type, which
are variously called \textit{graph monoids},
\textit{trace monoids} and \textit{free partially  commutative monoids},
arise naturally in the theory of computation, where they form a natural
model of concurrent processing \cite{Diekert97,Diekert95}.
Graph \textit{groups} can be used to model concurrent processing with
invertible operations; they also play an important role in combinatorial
group theory, where they are usually known as \textit{right-angled Artin
groups} \cite{Bestvina97,Brady01,Davis00,Howie99}.

The subgroup structure of graph groups has been extensively studied, with
extremely interesting results (see, for example, \cite{Bestvina97,Dicks99,Droms87}).
Likewise, there is considerable interest in submonoids of trace monoids. Of
particular importance in computer science are those submonoids of trace
monoids which are themselves trace monoids; an embedding of a trace monoid
into another is called a \textit{trace coding} \cite{Diekert97,Ochmanski88},
since it is the natural partially commutative analogue of a \textit{word
coding}. Trace codings have been extensively studied, with particular
attention paid to decidability questions \cite{Bruyere96,Chrobak87}.

Closely related to possible embeddings of graph groups or monoids, are the
possible commutation properties of subsets (and multisubsets) of groups and monoids.
These have been studied, in the group case, by Duncan,
Kazatchkov and Remeslennikov \cite{Duncan05}. Motivated by considerations
from algebraic geometry over groups, they associated to each finite graph
$\Gamma$ the class of groups which admit elements whose commutation
properties are described by the graph $\Gamma$.

In this paper, we consider the commutation properties of subsets of both
graph groups and graph monoids. In particular, we study certain key
graphs $\Gamma$ which have the property that a graph group
$G(\Omega)$ or graph monoid $M(\Omega)$ admits a subset whose commutation
properties are described by $\Gamma$ only when $\Omega$ contains an
embedded copy of $\Gamma$.
As a consequence, we obtain some negative results regarding embeddings
of both graph groups and graph monoids.

In addition to this introduction, this paper comprises four sections.
We begin, in Section~\ref{sec_basics}, by briefly introducing graph
groups and monoids, along with the notation and foundational results
which we shall need in the following sections.

Section~\ref{sec_groups} is devoted to graph groups. We show that a graph
group admits a subset whose commutation graph is a four-cycle if and only
if it contains an embedded (without chords) four-cycle. A consequence is
that a graph group contains a direct product of non-abelian free groups as
a subgroup only when its graph contains a four-cycle. This fact, which
was conjectured by Batty and Goda \cite{BattyGoda}, is of particular interest
because of the many properties which are shared by free groups and free
abelian groups but not by direct products of free groups. One
example is decidability of the algorithmic \textit{subgroup membership
problem}. A construction of Mikhailova \cite{Lyndon77,Mikhailova58} shows
that this problem is undecidable for direct products of free groups, and
hence for any graph group $G(\Gamma)$ where $\Gamma$ contains a chord-free
four-cycle. On the other hand, a recent result of
Kapovich, Weidmann and Myasnikov \cite{Kapovich05b} shows that the problem
is decidable for $G(\Gamma)$ when $\Gamma$ contain no chord-free cycles.
There remains the case of groups $G(\Gamma)$ where $\Gamma$ contains
chord-free cycles but not of length four. Our result shows that Mikhailova's
construction does not present an obstruction to decidability of subgroup
membership in these groups.

In Section~\ref{sec_monoids} we turn our attention to graph monoids. We
show that for certain graphs $\Gamma$, a graph monoid admits a subset with
commutation graph $\Gamma$ if and only if its graph contains an embedded
copy of $\Gamma$. As a consequence, we deduce a related restriction on
embeddings of direct products of free monoids.

In Section~\ref{sec_others} we ask what other graph groups and monoids
have similar properties. It transpires that the monoid results from
Section~\ref{sec_monoids} are best
possible, in the sense that every graph monoid which is \textit{not} a
direct product of free monoids of rank $1$ and $2$ admits an embedding
into a graph monoid without a corresponding embedding of graphs. This
contrasts with the group case, where it follows from a result of Droms,
Servatius and Servatius \cite{Droms92} that the graph group on the
three-edge line does not embed into a graph group without a corresponding
embedding of graphs. Finally, we give a combinatorial
construction
which embeds any member of a large class of graph groups into another graph
group, without a corresponding embedding of graphs.

\section{Graphs, Monoids and Groups}\label{sec_basics}

In this section, we briefly introduce the concepts, notation and
foundational results which will be required in the sections that
follow. We concentrate here on such of the theory as is common to
the monoid and group cases; ideas which are particular to groups or
monoids will be introduced in Sections~\ref{sec_groups} and
\ref{sec_monoids} respectively.

\subsection{Graphs}

By a \textit{graph} $\Gamma$ we mean a mean a set $V(\Gamma)$ of vertices
together with a reflexive, symmetric relation
$E(\Gamma) \subseteq V(\Gamma) \times V(\Gamma)$.
Two vertices $u, v \in V(\Gamma)$ are \textit{adjacent} if
$(u,v) \in E(\Gamma)$. The \textit{degree} $|\Gamma|$ of $\Gamma$ is
the cardinality of $V(\Gamma)$. The \textit{degree} of a vertex
$v$ in $\Gamma$, denoted $|v|_\Gamma$, is the number of vertices
adjacent to and distinct from $v$.

A \textit{morphism} from
a graph $\Gamma$ to a graph $\Omega$ is a map from the vertex set of
$\Gamma$ to that of $\Omega$ which preserves adjacency (but not
in general non-adjacency).
An \textit{embedding} of graphs is a morphism which
is injective on vertices and which preserves non-adjacency.
If $S$ is a subset of $V(\Gamma)$ then the \textit{subgraph of $\Gamma$
induced by $S$} is the graph with vertex set $S$ and edge set
$E(\Gamma) \cap (S \times S)$.

A
\textit{(connected) component} of a graph is a maximal set of vertices
such that every pair of vertices contained is connected by a path. A graph
is \textit{connected} if it has only one component, and
\textit{disconnected} otherwise.

Let $\Gamma$ and $\Omega$ be graphs with disjoint vertex sets. Then
the \textit{connected product} $\Gamma \times \Omega$ is the graph
with vertex set $V(\Gamma \times \Omega) = V(\Gamma) \cup V(\Omega)$ and
edge set
$$E(\Gamma \times \Omega) \ = \ E(\Gamma) \cup E(\Omega) \cup \left( V(\Gamma) \times V(\Omega) \right) \cup \left( V(\Omega) \times V(\Gamma) \right).$$

The
\textit{complement} $\ol{\Gamma}$ of $\Gamma$ is the graph with the same
vertex set as $\Gamma$, and in which two distinct vertices are adjacent
exactly if they are not adjacent in $\Gamma$.
A \textit{(co-connected) co-component} of
$\Gamma$ is a component of $\ol{\Gamma}$; the graph $\Gamma$ is called
\textit{co-connected} or \textit{co-disconnected} if $\ol{\Gamma}$ is
connected or disconnected respectively;

\begin{figure}\label{fig_graphs}
\begin{picture}(80,40)
\thicklines
\setloopdiam{10}
\Large
\setvertexdiam{1}

\letvertex A=(10,5)    \drawstate(A){}
\letvertex B=(10,15)    \drawstate(B){}
\letvertex C=(10,25)    \drawstate(C){}
\letvertex D=(10,35)    \drawstate(D){}
\drawundirectededge(A,B){}
\drawundirectededge(B,C){}
\drawundirectededge(C,D){}

\letvertex E=(30,15)    \drawstate(E){}
\letvertex F=(40,15)    \drawstate(F){}
\letvertex G=(40,25)    \drawstate(G){}
\letvertex H=(30,25)    \drawstate(H){}
\drawundirectededge(E,F){}
\drawundirectededge(F,G){}
\drawundirectededge(G,H){}
\drawundirectededge(H,E){}

\letvertex I=(60,15)    \drawstate(I){}
\letvertex J=(70,15)    \drawstate(J){}
\letvertex K=(60,25)    \drawstate(K){}
\letvertex L=(70,25)    \drawstate(L){}
\drawundirectededge(I,K){}
\drawundirectededge(J,L){}

\end{picture}
\caption{The graphs $L_3 = \ol{L_3}$, $C_4$ and $\ol{C_4} = E_{0,2}$.}
\end{figure}
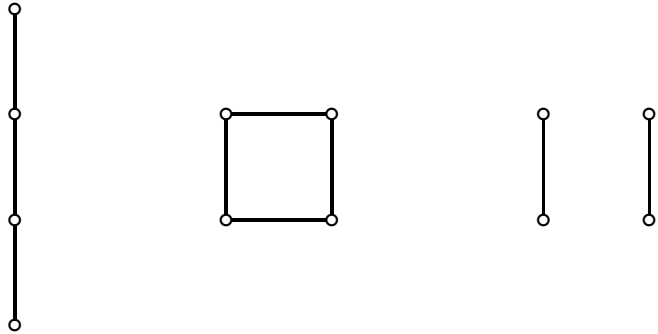

Figure~1 shows three examples of graphs which are important
in the study of
graph monoids and groups; they are the three-edge line $L_3$, the four-cycle
(or ``square'') $C_4$ and its complement $\ol{C_4}$. For clarity, we draw the graphs without the loops
at the vertices.
Note that $\ol{L_3}$ is isomorphic to $L_3$ --- we say that that $L_3$ is
\textit{self-complementary}. The complement graph $\ol{C_4}$ has
two connected components, each of which consists of two vertices joined by
an edge. More generally, for $i, j \geq 0$ we denote by $E_{i,j}$ the unique
graph with $i$ vertices of degree $0$ and $2j$ vertices of degree $1$, so
that $\ol{C_4} = E_{0,2}$. Another example which will be important for us
is the unique two-vertex disconnected graph $\ol{E_{0,1}} = E_{2,0}$.

\subsection{Graph Monoids and Groups}

Let $\Gamma$ be a graph. The \textit{graph monoid} $M(\Gamma)$ and \textit{graph group}
$G(\Gamma)$ are the monoid and group respectively defined by the
presentation
$$\langle \ V(\Gamma) \mid ab = ba \text{ for all } (a, b) \in E(\Gamma) \ \rangle.$$
There is an obvious embedding of $M(\Gamma)$ into $G(\Gamma)$, and it is
often convenient to regard the former as a submonoid of the latter.
It is also frequently useful to consider a set of monoid generators for
$G(\Gamma)$. With this in mind, we let
$$U(\Gamma) = \lbrace a, a^{-1} \mid a \in V(\Gamma) \rbrace$$
be a symmetrised set of generators for $G(\Gamma)$.
If $u \in U(\Gamma)^*$ then we denote by $\ol{u}$ the element of $G(\Gamma)$
(and hence, where appropriate, $M(\Gamma)$) represented.

Returning to our examples of graphs from above, it is easily seen that
$G(C_4)$ [respectively, $M(C_4)$] is isomorphic to a direct product two
free groups [monoids] of rank $2$. More generally, $G(\ol{E_{i,j}})$
[$M(\ol{E_{i,j}})$] is a direct product of $i$ free groups [monoids] of
rank $1$ and $j$ free groups [monoids] of
rank $2$. On the other hand, $G(L_3)$ [$M(L_3)$] is both freely and directly
indecomposable, and hence cannot be built up from free groups [monoids]
using only the operations of free and direct product. In fact, $L_3$ is
known to be the \textit{minimum} graph with the latter property, not only
in terms of number of vertices, but also with respect to embedding
\cite{Droms92}.

The \textit{length} of an element $g \in G(\Gamma)$ is the minimum length
of a word in the $U(\Gamma)^*$ representing $g$; a word of this length which
represents $g$ is called a \textit{reduced word} for $g$. A \textit{reduced
factorisation} for $g$ is an expression $g = g_1 \dots g_n$
where the sum length of the $g_i$s equals the length of $g$.
A \textit{prefix} [\textit{suffix}] of $g$ is an element $h \in G(\Gamma)$
which is the first factor [last factor] in some reduced factorisation for
$g$.

The \textit{support} of $g$ is the set of all vertices $t \in V(\Gamma)$
such that either $t$ or $t^{-1}$ or both occur in any (and hence in every)
reduced word for $g$. We write $t \in g$ to denote that $t \in V(\Gamma)$
lies in the support of $g \in G(\Gamma)$.
We say that two elements $u, v \in G(\Gamma)$ \textit{commute
totally} if every generator in the support of $u$ commutes with every
generator in the  support of $v$.
A generator $a \in U(\Gamma)$ is called \textit{central} if
$a$ has degree $|\Gamma| - 1$, that is, if $a$ represents a central
element in $G(\Gamma)$.

All of the definitions above apply by restriction to elements of
$M(\Gamma)$.

\subsection{Commutation Graphs}

Given a subset $S$ of a group $G$, the \textit{commutation graph of $S$
(in $G$)} is the graph with vertex set $S$, and an edge joining two
vertices exactly exactly if they commute in $G$.

Let $G$ be a group and $\Omega$ a graph. Following \cite{Duncan05},
we say that a group $G$ \textit{satisfies $\phi(\Omega)$} if there
exists a function $\sigma : V(\Omega) \to G$ with the property that
$u, v \in V(\Omega)$ are adjacent if and only if $\sigma(u)$ and
$\sigma(v)$ commute. Since $\sigma$ need not be injective, this is
in general slightly weaker than saying that $G$ has a subset with commutation
graph $\Omega$. However, with $\Omega$ countable and $G$ torsion-free
the two notions are easily seen to coincide, and it is this case
which will be of interest to us.

\section{Graph Groups and $C_4$}\label{sec_groups}

Our main aim in this section is to show that a graph group $G(\Gamma)$
does not satisfy $\phi(C_4)$ unless $\Gamma$ contains an embedded copy of
$C_4$. In particular, it follows $G(\Gamma)$ admits a subgroup isomorphic
to a direct product of free groups if and only if $\Gamma$ contains an
embedded copy of $C_4$; this proves part of a conjecture of Batty and
Goda \cite{BattyGoda}. In Section~\ref{sec_monoids}, we shall prove an
even stronger result in the monoid case.

Our main proof makes use of a theorem of Servatius \cite{Servatius89},
characterising centralizers of elements in graph groups; we begin by
briefly recalling some terminology and results from his paper. An element
$e \in G(\Gamma)$ is called \textit{cyclicly reduced} if it is of
minimal length amongst elements in its conjugacy class. It is easily
seen that every element of $G(\Gamma)$ can be written uniquely as a
reduced product $g = p h p^{-1}$ where $h$ is
cyclicly reduced. The element $h$ is called the
\textit{cyclic reduction} of $e$.

Now suppose $h$ is cyclicly reduced; and let $\Omega$ be the
subgraph of $\Gamma$ induced by the support of $h$. It is straightforward
to show that we can write $h$ in the form $h_1^{i_1} h_2^{i_2} \dots h_n^{i_n}$
where each $i_j$ is positive, each $h_j$ has support contained in different
co-connected component of $\Omega$, and no $h_j$ is a proper power. Moreover,
this expression is unique up reordering of the factors. The elements $h_i$
and their inverses are called the \textit{pure factors} of $h$. Notice
that the pure factors commute with one another.

\begin{theorem}[(The Centralizer Theorem, Servatius 1989 \cite{Servatius89})]\label{thm_centralizer}
Suppose $g = php^{-1}$ reduced, with $h$ cyclicly reduced. Then
$k$ commutes with $g$ if and only if $k$ can be written as
$p k_1 k_2 p^{-1}$ where $k_1$ is a product of pure factors of $h$,
and $k_2$ commutes totally with $h$.
\end{theorem}

The following combinatorial observation is probably well-known.
\begin{proposition}\label{prop_multiplyreduce}
Let $u, v \in G(\Gamma)$. Then there exist
reduced
factorisations $u = u'x$ and $v = x^{-1} v'$ such $u'v'$ is a reduced
factorisation for
$uv$.
\end{proposition}

\begin{proof}
Suppose false for a contradiction, and let $u$
and $v$ be counterexamples of minimal total length. Certainly $uv$ is
not a reduced product, or setting $u' = u$, $v' = v$ and $x = 1$ would
give the required properties.

Now let $\tilde{u}$ and $\tilde{v}$ be reduced words representing $u$ and
$v$ respectively. Since $uv$ is not reduced, we can write
$\tilde{u} \tilde{v} = \tilde{a} t \tilde{b} t^{-1} \tilde{c}$ for some
(positive or negative) generator $t$ and words $\tilde{a}$, $\tilde{b}$ and
$\tilde{c}$ such that $t$ commutes totally with $\tilde{b}$.

The factor
$t\tilde{b}t^{-1}$ cannot be contained entirely in $\tilde{u}$ or $\tilde{v}$,
since these are reduced words. It follows that the initial $t$ must
lie in $\tilde{u}$, and commute with every letter which comes after it
in $\tilde{u}$. By commuting it to the end, we may assume without loss of
generality that $t$ is the last letter in $\tilde{u}$. By a symmetrical
argument, we may assume also that $t^{-1}$ is the first letter in $\tilde{v}$.

Write $\tilde{u} = \tilde{a} t$ and $\tilde{v} = t^{-1} \tilde{b}$, and
let $a$ and $b$ be the elements represented by $\tilde{a}$ and $\tilde{b}$
respectively. Now
by the minimality
assumption, there are reduced factorisations $a = u' y$ and $b = y^{-1} v'$
such that $u'v'$ is a reduced factorisation for $ab = uv$. Now set
$x = yt$ to give $u = u'x$ and $v = x^{-1} v'$ as required.
\end{proof}

We shall need a number of other preliminary results. The first two are of a
rather technical nature.

\begin{proposition}\label{prop_onesided}
Suppose $y, t \in V(\Gamma)$ are two non-commuting vertices. Suppose
$t \notin uwv$ where $y \notin u$ and $y \notin v$ but $w$ is represented
by a reduced word beginning
and ending with a positive or negative occurrence of $y$. Then
$t \notin u$ and
$t \notin v$.
\end{proposition}
\begin{proof}
Suppose false for a contradiction, and let $u$ and $v$ be elements of
minimal total length such that the proposition fails, that is, such
that $t \in u$ or $t \in v$. By left-right symmetry, we can assume
without loss of generality that $t \in u$.

Let $\tilde{u}$, $\tilde{w}$ and $\tilde{v}$ be reduced words for
$u$, $w$ and $v$ respectively, where $\tilde{w}$ begins and ends with a
positive or negative occurrence of $y$.
Certainly $\tilde{u} \tilde{w} \tilde{v}$ is reducible, or we would have
$t \in uwv$, giving the
required contradiction. Hence, there must exist a factorisation
$\tilde{u} \tilde{w} \tilde{v} = \tilde{a} x \tilde{b} x^{-1} \tilde{c}$
where $\tilde{b}$ represents an element which commutes totally with the
positive or negative generator $x$.
Since the words $\tilde{u}$, $\tilde{v}$ and $\tilde{w}$ are reduced, the
factor $x\tilde{b}x^{-1}$ cannot lie wholly in any one of those words. Thus,
this factor must contain one end of $\tilde{w}$, and hence must contain $y$
or $y^{-1}$. Since $y \neq t$ and $y$ does not commute with $t$, it follows
that $x \neq t$ and $x \neq t^{-1}$. Since at least one occurrence of $x$
must lie in $\tilde{u}$ or $\tilde{v}$, we know also that $x \neq y$ and
$x \neq y^{-1}$

Now we can write $\tilde{a} \tilde{b} \tilde{c} = u' w' v'$ where $u'$, $w'$
and $v'$ are reduced
scattered subwords of $\tilde{u}$, $\tilde{w}$ and $\tilde{v}$ obtained by
deleting only occurrences of $x$ and $x^{-1}$, and the combined length of $u'$
and $v'$ is strictly less than that of $u$ and $v$. Moreover, it is clear
that $w'$ still begins and ends with a positive or negative occurrence of $y$.
Hence, by the minimality assumption, it
follows that $t \notin u'$ and $t \notin v'$,
and hence that $t \notin \tilde{u}$ and $t \notin \tilde{v}$. Thus,
$t \notin u$ and $t \notin v$, as required.
\end{proof}

\begin{lemma}\label{lem_anothertry}
Let $p \in G(\Gamma)$ and $t \in V(\Gamma)$, and suppose $r \in G(\Gamma)$
is of minimal length such that $p$ has a reduced factorisation of the form
$qtr$ or $qt^{-1}r$. Suppose $c \in G(\Gamma)$ is such that $t \notin p c p^{-1}$ and
the support of $c$ contains a generator which does not commute with $t$, and does not occur
in the support of $p$. Then $c$ has a reduced factorisation of the form
$c = r^{-1} t^{-1} d t r$ (if $p = qtr$) or $c = r^{-1} t d t^{-1} r$ (if $p = q t^{-1} r$).
\end{lemma}
\begin{proof}
We treat the case in which $p = qtr$; an entirely similar argument applies
when $p = qt^{-1}r$.

Let $y \in c$ be a generator which does not commute with $t$ and does not
occur in $p$. Let $\tilde{c}$ be a reduced word for $c$, and write
$\tilde{c} = \tilde{u} \tilde{w} \tilde{v}$ where $y \notin \tilde{u}$,
$y \notin \tilde{v}$ but $\tilde{w}$ begins and ends with a positive
or negative occurrence of $y$. Let $u$, $w$ and $v$ be the elements
represented by $\tilde{u}$, $\tilde{w}$ and $\tilde{v}$ respectively.

Now $t \notin pcp^{-1} = puwvp^{-1}$, so applying
Proposition~\ref{prop_onesided} we see that $t \notin pu = qtru$ and
$t \notin vp^{-1} = vr^{-1} t^{-1} q^{-1}$.

Now by Proposition~\ref{prop_multiplyreduce}, there exist reduced
factorisation s
$p = p'x$ and $u = x^{-1} u'$ such that $p'u'$ is a reduced factorisation
for $pu$. Now $t \notin pu = p'u'$, so clearly $p'$ does not
contain $t$. But $p'x = p = qtr$, so it follows easily from the minimality
assumption on $r$
that $x$ has a suffix $tr$, and hence that
$u$ has a prefix $r^{-1} t^{-1}$ as required. A symmetrical argument
shows that $v$ has a suffix $tr$.
\end{proof}

\begin{lemma}\label{lem_cyclicreducedimptotal}
Suppose $\Gamma$
does not contain an induced copy of $C_4$. Suppose further that $G(\Gamma)$
has a subset $S$ with commutation graph isomorphic to $C_4$, one of whose elements $a$ is
cyclicly reduced. Then $a$ commutes \textbf{totally} with
itself and with those other members of $S$ with which it commutes.
\end{lemma}
\begin{proof}
Suppose $\lbrace a, b, c, d \rbrace \subseteq G(\Gamma)$ has commutation graph isomorphic to
$C_4$, where $a$ fails to commute
with $c$, and $b$ fails to commute with $d$. Suppose further that $a$
is cyclicly reduced. Let $a_1, \dots, a_n$ be the pure factors of $a$.
Then by Theorem~\ref{thm_centralizer}, we can write
$b = b_1 b_2$ and $d = d_1 d_2$ where $b_1$ and $d_1$ are products of
pure factors of $a$, and $b_2$ and $d_2$ commute totally with $a$.

Notice that $b_1$ and $d_1$ commute, and both $b_1$ and $d_1$ commute
totally with both $b_2$ and $d_2$. Now if $b_2$ commuted with
$d_2$ then $b$ would commute with $d$, giving a contradiction. Hence,
$b_2$ does not commute with $d_2$.
It follows that  some vertex $x \in b_2$ fails to commute with some vertex
in $y \in d_2$.
Now $x$ and $y$ commute with every vertex in the support of $a$, so if two
vertices in the support of $a$ failed to commute then we would obtain a
four-cycle in $\Gamma$, giving a contradiction. Thus, $a$ must commute
totally with itself.

Now since the support of $b_1$ is contained in that of $a$, $b_1$ commutes
totally with $a$. We already know that $b_2$ commutes totally with $a$, so
it follows that $b$ commutes totally with $a$. By symmetry of assumption,
$d$ also commutes totally with $a$, as required.
\end{proof}

\begin{lemma}\label{lem_twocornercyclicreduced}
Suppose $\Gamma$ does not have an induced subgraph isomorphic to the
$C_4$ but $G(\Gamma)$ does satisfy $\phi(C_4)$. Then $\Gamma$ has a subset $S$
with commutation graph isomorphic to $C_4$, in which two commuting
elements are cyclicly reduced.
\end{lemma}
\begin{proof}
Suppose $\lbrace a, b, c, d \rbrace \in G(\Gamma)$ has commutation
graph isomorphic to $C_4$, where $a$ fails to commute
with $c$, and $b$ fails to commute with $d$. Clearly, by conjugating
the entire set, we may assume that one of these elements, say  $a$, is
cyclicly reduced. Now by Lemma~\ref{lem_cyclicreducedimptotal}, $b$
commutes totally with $a$.

Now suppose $b = p^{-1} b' p$ is a reduced factorisation where $b'$ is
cyclicly reduced. Clearly, the set
$\lbrace pap^{-1}, pbp^{-1}, pcp^{-1}, pdp^{-1} \rbrace$ also has
commutation graph isomorphic to $C_4$. But the support of $p$ is contained in that of
$b$, and hence $p$
commutes totally with $a$. It follows that $pap^{-1} = a$ and $pbp^{-1} = b'$
are both cyclicly reduced as required.
\end{proof}

\begin{theorem}\label{thm_groupmain}
Let $\Gamma$ be a graph containing no induced copy of $C_4$. Then
$G(\Gamma)$ does not satisfy $\phi(C_4)$.
\end{theorem}
\begin{proof}
Suppose for a contradiction that $G(\Gamma)$ satisfies $\phi(C_4)$.
By Lemma~\ref{lem_twocornercyclicreduced}, there is a subset
$\lbrace a, b, c, d \rbrace \subseteq G(\Gamma)$ with commutation graph
isomorphic to $C_4$ with $a$ and $b$ commuting and cyclicly
reduced. By Lemma~\ref{lem_cyclicreducedimptotal}, we may assume
without loss of generality that $a$ and $b$ commute totally with
themselves, with each other, and with $d$ and $c$ respectively.

Suppose $d = p^{-1} e p$ reduced with $e$ cyclicly reduced.
Then $c$ commutes with $d$ so by Theorem~\ref{thm_centralizer}, we can write
$c = p^{-1} c_1 c_2 p$ where $c_1$ is a product of
pure factors of $e$, and $c_2$ commutes totally with $e$.

Since $c$ does not commute with $a$, there must exist a letter
$x \in a$ which fails to commute with a letter $y \in c$.
Now $y$ must be in the support of at least one of $c_1$, $c_2$ and $p$. If
$y \in c_1$ or $y \in p$ then $y \in d$; but $d$ commutes totally
with $a$, so this contradicts the assumption that $x$ and $y$ do not
commute. Thus, we must have $y \in c_2$.

Since $b$ and $d$ do not commute, there are non-commuting vertices $s \in
b$ and $t \in d$. We know that $s$ and $t$ both commute with $x$, and that
$s$ commutes with $y$. We know also that $s$ and $t$ do not commute, and
that $x$ and $y$ do not commute. Since the graph is assumed to contain
no induced copy of $C_4$, it must be that $y$ does not commute with $t$.
Since $y \in c_2$ and $c_2$ commutes totally with $e$, we have $t \notin e$.
But $t \in d = pep^{-1}$, so we must have $t \in p$. However, since $t$ does
not commute with $s \in b$, and $c$ commutes totally with $b$, we must have
$t \notin c = p c_1 c_2 p^{-1}$.

Choose a reduced factorisation $qtr$ or $qt^{-1}r$ for $p$ such that
$r$ has minimal length. We have already observed that $y \in c_1 c_2$ does
not commute with $t$, and certainly $y \notin p$ or we would have $y \in d$
and $y$ would have to commute with $x$. Applying Lemma~\ref{lem_anothertry},
we see that $c_1 c_2$ has a
reduced factorisation of the form $r^{-1} t^{-1} f t r$. In particular,
the support of $t r$ is contained in the support of $c_1 c_2$.

Let $\Omega$ be the subgraph of $\Gamma$ induced by the support of $c_1 c_2$.
Notice that every vertex in the support of $c_1$ is connected to every vertex
in the support of $c_2$.
Thus, the support of $c_1$ and the support of $c_2$ are unions of disjoint
sets of co-connected components of $\Omega$.

Now we claim that the support of $t r$ is co-connected in $\Omega$. Indeed,
if not, then $r$ would contain letters from a co-component not containing
$t$; it would follow that we could commute these letters back through $t$,
contradicting the assumption that $r$ is of minimal length. It follows
that the support of $tr$ lies in a single co-connected component of
$\Omega$. In particular, the support of $tr$ is contained either in the
support of $c_1$, or in the support of $c_2$.

But $t$ cannot be in the support of $c_1$, since $c_1$ is a product of
pure factors of $e$, and $t$ is not in the support of $e$.
On the other hand, the final letter of $tr$ cannot be in the support of
$c_2$, since then it would commute with every letter in $e$, contradicting
the assumption that $pep^{-1} = qtrer^{-1}t^{-1}q^{-1}$ is a reduced
factorisation. This completes the proof.
\end{proof}
As an immediate corollary, we obtain a restriction on the graph groups which
contain a subgroup or submonoid isomorphic to a direct product of
non-abelian groups (or monoids).
\begin{corollary}\label{cor_groupmain}
Let $\Gamma$ be a graph not containing an embedded copy of $C_4$. Then
$G(\Gamma)$ has no subgroup [submonoid] isomorphic to direct product of
$2$ or more non-abelian groups [monoids].
\end{corollary}
\begin{corollary}\label{cor_groupembed}
If $\Gamma$ does not contain an embedded copy of $C_4$, then $G(\Gamma)$
does not contain an embedded copy of $G(C_4)$.
\end{corollary}

\section{Graph Monoids and $\ol{E_{i,j}}$}\label{sec_monoids}

In this section, we show that if $\Gamma$ is a graph in which every
vertex has degree $|\Gamma|-2$, that is, a graph of the
form $\ol{E_{0,j}}$, then a graph monoid $M(\Omega)$ satisfies
$\phi(\Gamma)$ only when $\Gamma$ embeds in $\Omega$. We deduce
also that a direct product of (abelian and non-abelian) free
monoids does not embed into a graph monoid without a
corresponding embedding of graphs.

We recall some standard definitions from the theory of graph monoids.
Given a graph $\Gamma$, we define a number of morphisms from $M(\Gamma)$
to free monoids of rank 1 and 2.
For each vertex $x \in V(\Gamma)$, let
$$\rho_x : M(\Gamma) \to \lbrace x \rbrace^*$$
be the map which deletes all symbols other than $x$. For each pair of
non-adjacent vertices $x$ and $y$, let
$$\sigma_{xy} : M(\Gamma) \to \lbrace x, y \rbrace^*$$
be the map which deletes all symbols other than $x$ and $y$.
The following well-known proposition says that any two distinct elements
of $M(\Gamma)$ are distinguished by at least one of the above morphisms;
a proof can be found in \cite{Diekert97}.
\begin{proposition}\label{prop_sigma_distinguishes}
Let $u, v \in V(\Gamma)^*$ be words in the vertices of $\Gamma$, such that $\ol{u}$ and $\ol{v}$ are distinct elements
of $M(\Gamma)$. Then either there exists a generator $x \in V(\Gamma)$ such that $\rho_x(u) \neq \rho_x(v)$ or there exist
non-commuting generators $x, y \in V(\Gamma)$ such that $\sigma_{xy}(u) \neq \sigma_{xy}(v)$.
\end{proposition}

This result, while elementary, is a key tool in the theory of graph 
monoids, and it will be central to our proofs below. We note that 
Proposition~\ref{prop_sigma_distinguishes} does not hold in the group 
case, with the obvious definitions of $\rho_x$ and $\sigma_{xy}$ as 
morphisms onto free groups of rank 1 and 2. For example, consider the 
graph $E_{1,1}$ with vertices $x$ of degree 0 and $y$ and $z$ of degree 1. 
Then the word $x y x^{-1} z x y^{-1} x^{-1} z^{-1} \in U(E_{1,1})^*$ does 
not represent the identity in $G(E_{1,1})$, but is not distinguished from 
the identity by a projection onto $1$ or $2$ generators. In fact, 
Proposition~\ref{prop_sigma_distinguishes} is a key reason why 
the theory of graph monoids is more straightforward than that of graph 
groups, and is 
why we obtain stronger results in the monoid case. One can formulate a 
related but more technical proposition concerning reduced words in the 
group generators \cite[Proposition~1]{Servatius89}, but this does not seem 
to be helpful for our purposes.

We proceed with a lemma characterising words which commute
in a graph monoid, in terms of the projections of the form $\sigma_{xy}$.
\begin{lemma}\label{lemma_sigma_and_commuting}
Let $u$ and $v$ be words in the vertices of $\Gamma$. Then $\ol{u}$ and $\ol{v}$ commute in $M(\Gamma)$ if
and only if for every pair of non-commuting vertices $x$ and $y$, there exists a word $w \in \lbrace x,y \rbrace^*$
which is not a proper power and integers $p, q \geq 0$ such that $\sigma_{xy}(u) = w^p$ and $\sigma_{xy}(v) = w^q$.
\end{lemma}
\begin{proof}
Certainly for any words $u$ and $v$ and vertex $x$ we have $\rho_x(uv) = \rho_x(vu)$,
so by Proposition~\ref{prop_sigma_distinguishes} we see that $\ol{u}$ and $\ol{v}$ commute if and only if for every pair of non-commuting vertices $x$ and
$y$ we have $\sigma_{xy}(uv) = \sigma_{xy}(vu)$.
Now $\sigma_{xy}$ is a morphism, so this is true if and only if for every $x$ and $y$,
$$\sigma_{xy}(u) \ \sigma_{xy}(v) \ = \ \sigma_{xy}(v) \ \sigma_{xy}(u),$$
that is, if $\sigma_{xy}(u)$ and $\sigma_{xy}(v)$ commute in the free monoid.
But clearly, this is the case if and only if $\sigma_{xy}(u)$ and $\sigma_{xy}(v)$ are powers
of a common subword, which can be chosen not to be a proper power.
\end{proof}

We need also the following lemma, which gives a necessary criterion for
distinct elements to commute.
\begin{lemma}\label{lemma_rho_distinguishes}
Let $u, v \in V(\Gamma)^*$ be words in the vertices of $\Gamma$, such that $\ol{u}$ and $\ol{v}$ are distinct elements
of $M(\Gamma)$. Suppose further that $\ol{u}$ and $\ol{v}$ commute in $M(\Gamma)$. Then there exists a
generator $x$ such that $\rho_x(u) \neq \rho_x(v)$.
\end{lemma}
\begin{proof}
Suppose not. Then by Proposition~\ref{prop_sigma_distinguishes}, there exist non-commuting generators $x$ and $y$
such that $\sigma_{xy}(u) \neq \sigma_{xy}(v)$.
But now by Lemma~\ref{lemma_sigma_and_commuting} there exists a word
$w \in \lbrace x,y \rbrace^*$ which is not a proper
power and integers $p, q \geq 0$ such that $\sigma_{xy}(u) = w^p$ and $\sigma_{xy}(v) = w^q$. Moreover,
since $\sigma_{xy}(u)$ and $\sigma_{xy}(v)$ are distinct, we must have $w$ non-empty and $p \neq q$.
Since $w$ is non-empty, it must contain either an $x$ or a $y$. Suppose without loss of generality that
it contains $k \geq 1$ occurrences of the letter $x$. Then $\sigma_{xy}(u)$ contains $kp$ occurrences
of $x$, while $\sigma_{xy}(v)$ contains $kq$ occurrences of $x$. It follows that $u$ and $v$ contain
$kp$ and $kq$ occurrences of $x$ respectively, so that $\rho_x(u) \neq \rho_x(v)$. This contradicts our
supposition and hence completes the proof.
\end{proof}

\begin{lemma}\label{lemma_commute_or_miss}
Suppose $u, v \in V(\Gamma)^*$ are words in the vertices of $\Gamma$, such that $\ol{u}$ and $\ol{v}$ commute.
Let $x \in V(\Gamma)$ be a vertex which occurs in $u$. Then either $x$ occurs in $v$, or $x$ commutes with
every letter which occurs in $v$.
\end{lemma}
\begin{proof}
Suppose for a contradiction that $x$ does not occur in $v$, and does not commute with some letter $y$
which occurs in $v$. By Lemma~\ref{lemma_sigma_and_commuting}, there exists a word $w \in F(x,y)$ and integers
$p, q \geq 0$ such that $\sigma_{xy}(u) = w^p$ and $\sigma_{xy}(v) = w^q$. Now $u$ contains the letter
$x$, so $\sigma_{xy}(u) = w^p$ contains the letter $x$, so $w$ must contain the letter $x$. On the other hand,
$\sigma_{xy}(v) = w^q$ does not contain the letter $x$, so we must have $q = 0$ and
$\sigma_{xy}(v) = \epsilon$. But $v$ contains an occurrence of $y$, so $\sigma_{xy}$ contains an
occurrence of $y$, and in particular is non-empty. This gives the required contradiction.
\end{proof}

\begin{lemma}\label{lem_monoidinduction}
Suppose $M(\Gamma)$ satisfies $\phi(\Omega \times \ol{E_{0,1}})$. Then
$\Gamma$ has an induced subgraph isomorphic to $\Gamma_1 \times \ol{E_{0,1}}$
where $\Gamma_1$ satisfies $\phi(\Omega)$.
Moreover, if $M(\Omega \times \ol{E_{0,1}})$ embeds into $M(\Gamma)$ then
$M(\Omega)$ embeds into $M(\Gamma_1)$.
\end{lemma}
\begin{proof}
Let $S$ be a subset of $M(\Gamma)$ with commutation graph isomorphic to
$\Omega \times \ol{E_{0,1}}$, let $e, e' \in S$ be the elements which
map to the vertices of $\ol{E_{0,1}}$ under this isomorphism, and let
$S' = S \setminus \lbrace e, e' \rbrace$. Thus, $e$ and $e'$
commute with every element of $S'$, but not with each other.

Let $\Gamma_1$ be the subgraph induced by the set of all vertices in
$V(\Gamma)$ which occur in
the support of elements in $S'$. It is immediate from the definition that
$\Gamma_1$ satisfies $\phi(\Omega)$.
Moreover, if $S$ actually generates a submonoid isomorphic to
$M(\Omega \times \ol{E_{0,1}})$ and $e$ and $e'$ are chosen appropriately, then
$M(\Omega)$ embeds in $M(\Gamma_1)$.

By Lemma~\ref{lemma_sigma_and_commuting} we may choose vertices
$y, z \in V(\Gamma)$ such that
$\sigma_{y z}(e)$ and $\sigma_{yz}(e')$ are
not powers of a common subword. Let $\Gamma_2$ be the subgraph induced
by the vertex set $\lbrace y, z \rbrace$.
To prove the lemma, it will suffice to show that $V(\Gamma_1)$ and $V(\Gamma_2)$ are
disjoint, and that every vertex in $\Gamma_2$ is adjacent to every vertex
in $\Gamma_1$.

First, we claim that $V(\Gamma_1)$ and $V(\Gamma_2)$ are disjoint, that is, that
$y$ and $z$ do not lie in $V(\Gamma_1)$. Indeed suppose for a contradiction
that at least one of
them does, and
let $d \in S'$ be an element with support containing $y$ or $z$.
Then $\sigma_{y z}(d)$ is non-empty.
We know that $d$ commutes with $e$ and $e'$, so by Lemma~\ref{lemma_sigma_and_commuting} we have
$\sigma_{y z}(e)$ and
$\sigma_{y z}(d)$ are powers of a common subword, and likewise
that $\sigma_{yz}(e')$ and $\sigma_{yz}(d)$.
But the relation of being powers of a common subword is transitive through
non-empty words, so
it would follow that $\sigma_{yz}(e)$ and $\sigma_{yz}(e')$
are powers of a
common subword, giving the required contradiction.

Now since every element in $S'$ commutes with $e$ and with $e'$,
Lemma~\ref{lemma_commute_or_miss} tells us that every letter in $\Gamma_1$
commutes with $y$ and with $z$. This completes the proof.
\end{proof}

An inductive applications of Lemma~\ref{lem_monoidinduction} leads to
the first main theorem of this section. We note that the case $j = 1$
can also be obtained as a consequence of Theorem~\ref{thm_groupmain}.
\begin{theorem}\label{thm_monoidmain}
Let $j \geq 0$. Then $M(\Gamma)$ satisfies $\phi(\ol{E_{0,j}})$ if and
only if $\ol{E_{0,j}}$ embeds in $\Gamma$.
\end{theorem}
\begin{proof}
Suppose the direct implication is false, and let $j \geq 0$ be minimal
such that there exists a graph $\Gamma$
such that $M(\Gamma)$ satisfies
$\phi(\ol{E_{0,j}})$ but $\ol{E_{0,j}}$
does not embed into $\Gamma$.
 Certainly $j \neq 0$, since $\ol{E_{0,0}}$ is
the graph with no vertices, which certainly embeds into $\Gamma$.

Otherwise, we have $\ol{E_{0,j}} = \ol{E_{0,j-1}} \times \ol{E_{0,1}}$, so
by Lemma~\ref{lem_monoidinduction}, we see that $\Gamma$ has a subgraph
isomorphic to $\Gamma_1 \times \ol{E_{0,1}}$
where $\Gamma_1$ satisfies $\phi(\ol{E_{0,j}})$.
By the minimality assumption on $j$, $\ol{E_{0,j-1}}$ embeds into $\Gamma_1$,
and it follows that $\ol{E_{0,j-1}} \times \ol{E_{0,1}} = \ol{E_{0,j}}$ embeds
into $\Gamma_1 \times \ol{E_{0,1}}$, and hence into $\Gamma$ as required.

The converse implication is immediate.
\end{proof}

Before proving our second main theorem of this section, we need the following
preliminary step. We remark that Batty and Goda \cite{BattyGoda}
have observed that an analogous result holds for groups.
\begin{proposition}\label{prop_commutative_unconcealable}
Let $\Gamma$ be a graph not containing an induced subgraph isomorphic
to the complete graph $\ol{E_{n,0}}$ on $n$ vertices.  Then $M(\Gamma)$
does not have a submonoid
isomorphic to the free commutative monoid of rank $n$.
\end{proposition}
\begin{proof}
Suppose false for a contradiction, and let $\Gamma$ be a graph of
minimal degree such that the claim fails. 
 Let
$u_1, \dots u_n \in V(\Gamma)^*$ be words in the vertices of $\Gamma$ such
that the corresponding elements $\ol{u_1}, \dots, \ol{u_n} \in M(\Gamma)$
generate a free commutative monoid $N$ of rank $n$.

It follows from elementary linear algebra that the free commutative
monoid of rank $n$ does not
embed into a free commutative monoid of rank less than $n$, so we may assume
that $\Gamma$ is not a complete graph and choose non-adjacent
vertices $x, y \in V(\Gamma)$.

It follows from Lemma~\ref{lemma_sigma_and_commuting} that there exists a 
word $r \in \lbrace x, y \rbrace^*$ such that each $\sigma_{xy}(u_i)$ is of
the form $r^q$ 
for some $q \geq 0$. Since $\Gamma$ is of minimal degree, every vertex of 
$\Gamma$ must occur in some $u_i$. In particular, $x$ and $y$ must each
occur in some $u_i$ and so they must both occur in $r$.

We define a morphism $f : M(\Gamma) \to M(\Gamma)$ by letting $f(w)$
be obtained from $w$ by deleting all occurrences of the generator $x$.
We claim that this morphism is injective on $N$. Indeed, suppose
$w, w' \in V(\Omega)^*$ represent distinct elements of $N$. Then by
Lemma~\ref{lemma_rho_distinguishes}, we have $\rho_a(w) \neq \rho_a(w')$ for
some $a \in V(\Omega)$. We claim that we may assume without loss of generality that $a \neq x$. Indeed,
by our observations above, $\sigma_{xy}(w) = r^p$ and $\sigma_{xy}(w') = r^q$
for some $p, q \geq 0$.
If $\rho_x(w) \neq \rho_x(w')$ then we must have $p \neq q$. But since $r$
contains at least one
occurrence of $y$, it follows that $\rho_y(w) \neq \rho_y(w')$ so 
we can instead take $a = y$.

Now we have
$$\rho_a(f(w)) = \rho_a(w) \neq \rho_a(w') = \rho_a(f(w')),$$
so that $\ol{f(w)} \neq \ol{f(w')}$. This proves the claim that
$f$ is injective.

Now since the image $f(N)$ is contained within the induced subgraph with
vertex set $V(\Gamma) \setminus \lbrace x \rbrace$, this contradicts the
minimality assumption on $\Gamma$ and completes the proof.
\end{proof}

We are now ready to prove the following theorem.
\begin{theorem}\label{thm_monoidembedding}
Let $i, j \geq 0$. Then $M(\Gamma)$ has a submonoid isomorphic to a direct
product of $i$ rank $1$ free monoids and $j$ non-abelian free monoids if
and only if $\ol{E_{i,j}}$ embeds in $\Gamma$.
\end{theorem}
\begin{proof}
Suppose $M(\Gamma)$ has a submonoid isomorphic to a direct product of
product of $i$ rank $1$ free monoids and $j$ non-abelian free monoids. 
Then clearly, $M(\Gamma)$ has a submonoid isomorphic to $M(\ol{E_{i,j}})$.
Notice that
$$\ol{E_{i,j}} = \ol{E_{i,0}} \times \ol{E_{0,j}} = \ol{E_{i,0}} \times \ol{E_{0,1}} \times \dots \times \ol{E_{0,1}}.$$
By an inductive application of Lemma~\ref{lem_monoidinduction}, we
deduce that $\Gamma$ has a subgraph isomorphic to $\Gamma_1 \times \ol{E_{0,j}}$
where the free commutative monoid $M(\ol{E_{i,0}})$ of rank $i$ embeds in
$M(\Gamma_1)$.

Now by Proposition~\ref{prop_commutative_unconcealable} we deduce
that $\Gamma_1$ contains a complete subgraph with $i$ vertices. It
follows that $\Gamma_1 \times \ol{E_{0,j}}$ has a induced subgraph isomorphic
to $\ol{E_{i,j}}$, and hence so does $\Gamma$.
\end{proof}

\section{Other Graph Monoids and Groups}\label{sec_others}

It seems natural to ask whether similar results hold for other graphs,
that is, whether there are other graphs $\Gamma$ with the property that
$G(\Omega)$ or $M(\Omega)$ satisfies $\phi(\Gamma)$ only when $\Gamma$
embeds in $\Omega$.

A related, but weaker, property has been considered by
Batty and Goda \cite{BattyGoda}. They call a graph group $G(\Gamma)$
\textit{unconcealable} if it embeds into a graph group $G(\Omega)$ only
when $\Gamma$ embeds into $\Omega$. They observe that the free group of
rank $2$ and all free abelian groups are unconcealable, and conjecture
that direct products of free groups of rank $1$ and $2$ also have this
property.
Thus, our Corollary~\ref{cor_groupembed} proves one case of their conjecture;
the general case remains open.

The notion of unconcealability applies equally naturally to monoids, and
our Theorem~\ref{thm_monoidmain} is the natural monoid-theoretic analogue
of Batty and Goda's conjecture. In fact, in the monoid case, it transpires
that this result is best possible, in the sense that all graphs not covered
by that theorem admit concealments.
\begin{proposition}
If $M(\Gamma)$ is unconcealable then $\Gamma = \ol{E_{i,j}}$ for some
$i, j \geq 0$.
\end{proposition}
\begin{proof}
Suppose $M(\Gamma)$ is unconcealable, and consider the direct product of
the projections $\sigma_{xy}$ and $\rho_x$. This is an injective
(by Proposition~\ref{prop_sigma_distinguishes}) morphism
from $M(\Gamma)$ to a direct product of free monoids of rank $1$ and $2$,
that is, an embedding of $M(\Gamma)$ into a monoid of the form
$M(\ol{E_{i',j'}})$. Since $M(\Gamma)$ is unconcealable, $\Gamma$ must
embed into $\ol{E_{i',j'}}$. It now follows easily that $\Gamma$ is of
the form $\ol{E_{i,j}}$.
\end{proof}
In the group case, the lack of a counterpart to
Proposition~\ref{prop_sigma_distinguishes} once again means that things
are not so straightforward. In general, it is not clear exactly which graph
groups are unconcealable or have our stronger property.

Recall that an \textit{assembly} group is a graph group which can be built
up from copies of $\Z$ using free and direct products. Droms, Servatius
and Servatius \cite{Droms92} have shown that no non-assembly graph group
embeds into an assembly group. They observe also that $G(\Gamma)$ (with
$\Gamma$ finite) is an
assembly group if and only if $\Gamma$ contains
no embedded copy of $L_3$. Thus, their result can be
interpreted as saying that $G(L_3)$ is unconcealable. It seems natural
also to ask if this graph has our stronger property.
\begin{question}
Is there an assembly group satisfying $\phi(L_3)$?
\end{question}

The rest of this
section is devoted to a combinatorial construction which yields a
concealment for a large number of graph groups (and monoids).
Specifically, we show that for $G(\Gamma)$ to be unconcealable it is
necessary either that every vertex has degree $|\Gamma|-2$ or more
(that is, $\Gamma = \ol{E_{i,j}}$ for some $i, j \geq 0$) or that
$\Gamma$ has vertices of degree $|\Gamma|-2$ and $|\Gamma|-3$.

Let $\Gamma$ be a graph which does not satisfy this condition, that is,
which has a vertex of degree $|\Gamma|-3$ or
less, but does not have vertices of degree both $|\Gamma|-2$ and $|\Gamma|-3$.
Let $e$ be a vertex of maximal degree amongst those vertices having degree
$|\Gamma| - 3$ or less, and let $f$ and $g$ be vertices which are not adjacent
to $e$.

Let $e_0$ and $e_1$ be new symbols not in $V(\Gamma)$ and define a new graph $\Omega$
with
$$V(\Omega) \ = \ \left[ V(\Gamma) \setminus \lbrace e \rbrace \right] \cup \lbrace e_0, e_1 \rbrace, \text{ and}$$
\begin{align*}
E(\Omega) \ = \ &V(\Omega) \cup \lbrace (e_0, a), (a, e_0), (a, e_1), (e_1, a) \mid (e, a) \in E(\Gamma) \rbrace \\
          &\cup \lbrace (e_0, f), (f, e_0), (e_1, g), (g, e_1) \rbrace.
\end{align*}
We claim that $G(\Gamma)$ and $M(\Gamma)$ are concealed in
$G(\Omega)$ and $M(\Omega)$ respectively. We begin by showing
that $\Gamma$ is not an induced subgraph of $\Omega$.
\begin{proposition}\label{prop_notinduced}
$\Gamma$ does not embed in $\Omega$.
\end{proposition}
\begin{proof}
Suppose for a contradiction that $\Omega$ has an induced subgraph
$\Sigma$ which is isomorphic to $\Gamma$. Since $|\Gamma| = |\Omega| - 1$,
$\Sigma$ must be induced by deleting one vertex from $\Omega$; call this
vertex $v$.

By construction,
$\Omega$ has $|e|_\Gamma + 2$ more edges than $\Gamma$.
In order for $\Sigma$ to have the same number of edges as $\Gamma$, it must
be that $|v|_\Omega = |e|_\Gamma + 2$. In particular $v$
cannot be $e_0$ or $e_1$, both of which are constructed to have
degree $|e|_\Gamma + 1$ in $\Omega$. It follows that $v$ is a vertex from $\Gamma$.
Now by the construction of $\Omega$, $|v|_\Gamma$ must be either
$|v|_\Omega - 1 = |e|_\Gamma + 1$ (if $v = f$, $v = g$ or $(v, e) \in E(\Gamma)$)
or $|v|_\Omega = |e|_\Gamma + 2$ (otherwise).
Hence, by the maximality assumption on $|e|_\Gamma$, either
$|v|_\Gamma = |\Gamma| - 1$ or $|v|_\Gamma = |\Gamma| - 2$.

Suppose first that $|v|_\Gamma = |\Gamma| - 1$, that is, that $v$ is
central in $\Gamma$. Note that $v$ cannot be $f$ or $g$, since
neither commute with $e$ in $\Gamma$. Now suppose a vertex $x$ is
central in $\Sigma$. Certainly
$x \neq e_0$, since $e_0$
does not commute with $g$ in $\Omega$, and $g$ remains in the
induced subgraph. By a symmetrical argument, $x \neq e_1$, so $x$
must also be a vertex in $\Gamma$.
Moreover, $x$ commutes with every vertex in $\Sigma$ and also
with $v$. It follows easily that $x$ is central in
$\Gamma$. We have shown that every central vertex in $\Sigma$ is a central
vertex in $\Gamma$, and we know also that $v$ is a central vertex in
$\Gamma$. But now $\Sigma$ has strictly fewer central vertices than
$\Gamma$, which
contradicts the assumption that $\Sigma$ is isomorphic to $\Gamma$.

Now suppose that $|v|_\Gamma = |\Gamma| - 2$. We have already seen
that $|e|_\Gamma$ is either $|v|_\Gamma - 1$ or $|v|_\Gamma - 2$.
But by our original assumptions, $\Gamma$ cannot contain a vertex of
degree $|\Gamma|-3$, so it must be that $|e|_\Gamma = |\Gamma|-4$.
It follows from the construction of $\Omega$ that $\Omega$ has the same
number of vertices of degree greater than or equal to $|\Gamma| - 2$ that
$\Gamma$ does. But $v$ has degree $|\Gamma| - 2$ and is missing
from $\Sigma$. Hence, $\Sigma$ has
strictly fewer vertices of degree greater than or equal to
$|\Gamma| - 2$ than $\Gamma$, which again gives the required
contradiction.
\end{proof}
Now considering for example the subset
$(\Gamma \setminus \lbrace e \rbrace) \cup \lbrace e_0e_1e_0e_1 \rbrace$
it is clear that $M(\Omega)$ has a subset with commutation graph $\Gamma$, and
so both $M(\Omega)$ and $G(\Omega)$ satisfy $\phi(\Gamma)$.
In fact, we can go further. Define a monoid morphism
$$
\tau : U(\Gamma)^* \to U(\Omega)^*, \ \ 
\tau(x) = \begin{cases}
               e_0 e_1 e_0 e_1 &\text{ if } x = e; \text{ or} \\
               e_1^{-1} e_0^{-1} e_1^{-1} e_0^{-1} &\text{ if } x = e^{-1}; \text{ or} \\
               x               &\text{ if } x \in U(\Gamma) \setminus \lbrace e, e^{-1} \rbrace.
            \end{cases}
$$
It is immediate from the definition of $\Omega$ that $\tau$ respects the
defining relations in $G(\Gamma)$, and so induces a well-defined morphism
$$\ol{\tau} : G(\Gamma) \to G(\Omega), \ \   \ol{w} \mapsto \ol{\tau(w)}.$$
A straightforward but technical argument shows that $\ol{\tau}$ is injective,
thus completing the proof of the following.
\begin{theorem}
Let $\Gamma$ be a graph which has a vertex of degree $|\Gamma|-3$ or
less, but does not have both a vertex of degree $|\Gamma|-3$ and a vertex
of degree $|\Gamma|-2$. Then there exists a graph $\Omega$ with
$|\Omega| = |\Gamma| + 1$ such that $G(\Gamma)$ embeds into $G(\Omega)$,
but $\Gamma$ does not embed into $\Omega$.
\end{theorem}

\section*{Acknowledgements}

The research documented here was started while the author was at
Carleton University, supported by the Leverhulme Trust. It was completed,
and this paper written, at Universit\"at Kassel, where it was supported by
a Marie Curie Intra-European Fellowship within the 6th European
Community Framework Programme.
The author would like to thank Mike Batty and Keith Goda for many
helpful discussions. He would also like to thank Kirsty for all her
support and encouragement.

\bibliographystyle{plain}
\bibliography{mark}

\def\cprime{$'$} \def\cprime{$'$}

\end{document}